\documentclass{amsart}
\usepackage{graphicx}
\graphicspath{{figures/}}
\usepackage{amssymb}
\usepackage[font={footnotesize,it}]{caption}
\usepackage{comment}
\usepackage{url}
\usepackage{hyperref}
\usepackage{caption}
\usepackage{mathtools}
\usepackage[colorinlistoftodos,prependcaption,textsize=tiny]{todonotes}

\numberwithin{theorem}{section}
\theoremstyle{remark}

\DeclareMathOperator\Arg{Arg}

\begin{document}
\title[Computation of generalized Stieltjes constants]{High precision computation and a new asymptotic formula for the generalized Stieltjes constants}
\author{Sandeep Tyagi}
\thanks{}
\email{tyagi\_sandeep@yahoo.com}
\keywords{Stieltjes constants, Hurwitz zeta function, Riemann zeta function,  double exponential method}
\subjclass[2010]{Primary 11M06, 11Y16; Secondary 68Q25.}

\begin{abstract}
We provide an efficient method to evaluate the generalized Stieltjes constants $\gamma_n(a)$ numerically to arbitrary accuracy for large $n$ and $n \gg |a|$ values. The method uses an integral representation for the constants and evaluates the integral by applying the double exponential (DE) quadrature method near the saddle points of the integrands.  Further, we provide a highly accurate asymptotic formula for the generalized Stieltjes constants.
\end{abstract}

\maketitle

\section{Stieltjes Constants}
The generalised Stieltjes constants (SCs) are defined in terms of the Laurent expansion of the Hurwitz zeta function $\zeta(s,v)$ around $s=1$:
\begin{align}
\zeta(s,v) = \frac{1}{s-1} + \sum_{n=0}^{\infty} \frac{(-1)^n}{n!} \gamma_n(v) (s-1)^n \quad \quad s \ne 1,   
\label{zsc}
\end{align}
where $s$ and $v$ are complex numbers and $v$ is restricted to region $\text{Re}(v) \ge \frac{1}{2}$. 

There is a long history of scientific interest in these constants (see \cite{johansson2019computing} and references therein). They are important in number theory.  In particular, asymptotic bounds of $|\gamma_n(1)|$ for $n \gg 1$ have been used to identify the zero free region of the Riemann zeta function \cite{eddin2017applications}.

There are various numerical methods known to evaluate these constants \cite{johansson2015rigorous,  coffey2006new, coffey2014series,  knessl2011effective,  fekih2014new,  adell2017fast,  maslanka2022high} for small $n$ values.  However, an accurate estimation of these constants has been a challenge for large $n$ especially when the accuracy desired  for $\gamma_n(v)$ is greater than $\log_{10}(n)$. Significant progress in this direction was made by Johansson \& Blagouchine \cite{johansson2019computing} who used the following representation \cite{blagouchine2015theorem} to calculate $\gamma_n(v)$ for arbitrary $n$:
\begin{align}
\gamma_n(v) = -\frac{\pi}{2(n+1)} \int_0^{\infty} \frac{\log^{n+1}(v-\frac{1}{2}-i x)+\log^{n+1}(v-\frac{1}{2}+i x)}{\cosh^2(\pi x)} dx.
\label{sc}
\end{align}

They use \autoref{sc} to compute the SCs numerically using Gaussian quadrature in ball arithmetic. However, these calculations become slower as $n$ becomes large.  In this paper we improve on their method and show that $\gamma_n(v)$ can be calculated in a straightforward way by applying the double exponential (DE) method \cite{mori2001double, takahasi1974double, tyagi2022double} around the saddle point without the need to use the ball arithmetic.  

Asymptotic formulas for Stieltjes constants  have been given by Matsuoka \cite{matsuoka1985generalized},  Knessl \& Coffey \cite{knessl2011effective} and quite recently by Masalanka \cite{maslanka2022asymptotic}.   Asymptotic expansion for generalized SCs is given by Knessl \& Coffey \cite{knessl2011asymptotic} and Paris \cite{paris2015asymptotic}.  In this paper we derive a highly accurate asymptotic expansion for the generalized SC constants.

\section{Numerical Algorithm}
Following Johansson \& Blagouchine \cite{johansson2019computing}, we define $a=v-\frac{1}{2}$,
\begin{align}
f(x) \,\equiv\, \frac{\log^{n+1}(a+ix)}{\cosh^2 \! \pi x }\,, \qquad I_n(a) \,\equiv\, \int_0^{\infty} \!\! f(x) \, dx.
\label{eq:f}
\end{align}
The SCs are then expressed as
\begin{align}
\gamma_n(v) = 
-\frac{\pi}{(n+1)} \cdot
\begin{cases}
2 \Re(I_n(a)), & \Im(a) = 0  \\[2mm]
I_n(a) + \overline{I_n(\overline{a})}, & \Im(a) \ne 0.  
\end{cases}
\label{eq:gammafromi}
\end{align}
where ``$\,\overline{\phantom{m}}\,$'' stands for the complex conjugate. The integrand is written as
\begin{align}
f(x) = \exp\left(g(x)\right) h(x),
\label{eq:fgh1}
\end{align}
where
\begin{align}
g(x) = (n+1) \log\left(\log\left(a + ix\right)\right) - 2 \pi x, \quad h(x) = (1 + \tanh\pi x)^2.
\label{eq:fgh2}
\end{align}
Assuming that $n \gg |a|$, the function $\exp(g(x))$ has a single saddle point in the right half-plane and this can be found by solving the equation $g'(\omega) = 0$. This leads to
\begin{align}
(n+1) + 2 \pi i \left(a + i \omega \right) \log\left(a + i \omega \right) = 0,
\end{align}
which can be solved for $w$:
\begin{equation}
\omega = i \left( a - \frac{u}{W_0(u)}\right), \quad u = \frac{(n+1) i}{2 \pi},
\label{eq:owega}
\end{equation}
where $W_0(u)$ is the principal branch of the Lambert $W$ function \cite{corless1996lambertw} .

We apply the DE method to evaluate the integral $I_n(a)$ in \autoref{eq:f}. For integration, we use a contour passing through the saddle point roughly along the steepest direction of descent.  The direction of the contour, $\epsilon = \exp(i \phi)$,  at the saddle point can be obtained by finding the second derivative of $g(z)$ at $z=\omega$,
\begin{align}
g''(\omega)=(n+1) \frac{1+\log^{-1}(a+i \omega)}{(a+i \omega)^2 \log(a+i \omega)},
\label{eq:double}
\end{align}
and setting $\phi = (\pi- \theta)/2$ where  $\theta =\Arg{(\omega^2 g''(\omega))}$. Thus,  the contour of integration can be set as  $x(y)$:
\begin{align}
x(y)=\omega e^{1+\epsilon y-e^{-\epsilon y}}.
\end{align}
This choice ensures that saddle point is located at $y=0$ and that the range of integration changes from $(0,\infty)$ to $(-\infty, \infty)$.  In practice we set $\epsilon=1$ which works well. This is because numerically the real part $g''(\omega)$ dominates over its imaginary part.

\begin{figure}[ht]
\centering
\includegraphics[width=0.9\textwidth]{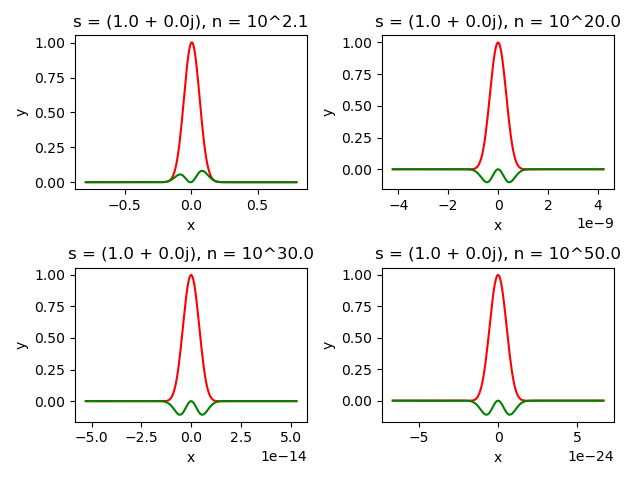}
\caption{Double exponential decay of the integrand around the saddle point for $v=1$ and various $n$ values. The real (red) and imaginary parts (green). are shown  Note the scaling factor over the x-axis showing a sharp decay away from the saddle point.}
\label{fig:j1}
\end{figure}

\begin{figure}[ht]
\centering
\includegraphics[width=0.9\textwidth]{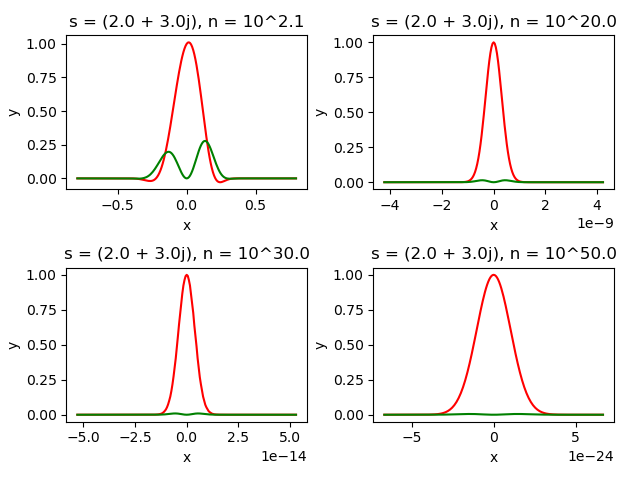}
\caption{Double exponential decay of the integrand around the saddle point for $v=2+3 i$ and various $n$ values. The real (red) and imaginary parts (green). are shown  Note the scaling factor over the x-axis showing a sharp decay away from the saddle point.}
\label{fig:j2}
\end{figure}

Note that the integrand shows a DE decay along the positive and negative $y$ directions. We show the profile of function $f(y)/f(0)$ near the saddle point, $y=0$,  in \autoref{fig:j1} for $v=1$ and different $n$ values.  A similar plot for $v=2+3 i$ is shown in  \autoref{fig:j2}. As expected, both the real and imaginary parts show a very smooth behaviour near the saddle point and sharp decay away from it. 

The integration region is cut off at $\pm q$. For a desired accuracy of $10^{-N}$, the range $q$ can be found roughly by $q=\sqrt{N \ln(10)/ |g''(\omega) \omega^2|}$, where the double derivative of $g(w)$ is given in \autoref{eq:double}. We discretize the integration region $(-q,+q)$ into $M=201$ parts.  This value of $M$ is sufficient for an accuracy of 100 digits for $n \ge 500$. 

\begin{figure}[ht]
\centering
\includegraphics[width=0.9\textwidth]{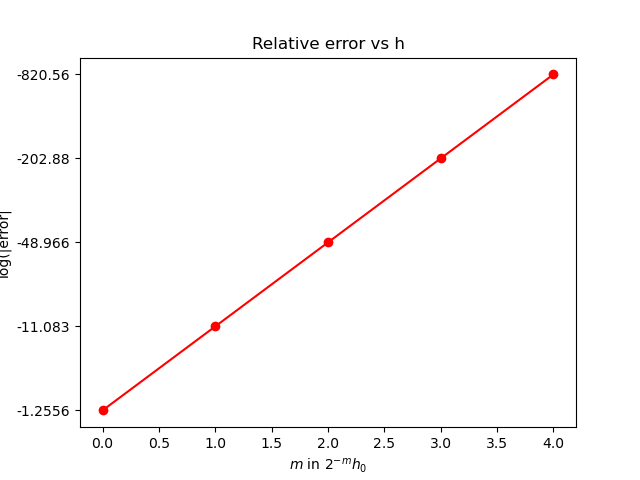}
\caption{Relative error as a function of discretization size $h$ is shown.  Here $h = h_0 2^m$ and $h_0=1.2146 \times 10^{-49}$ and integration cut off was at $q=6.067 \times 10^{-48}$.}
\label{fig:j3}
\end{figure}

Every halving of $h$ is guaranteed to at least double the accuracy achieved,  and a convergence in results can be detected easily.   We check the convergence of accuracy with respect to discretization size in \autoref{fig:j3}.  We work for $n=10^{100}$ case.  First we choose an appropriate $q$ value to ensure we can obtain an accuracy of more than 1000 digits.  We choose M=3201 which was found to be sufficient to obtain more than 1000 digits accuracy in this case.  We then take the result obtained for these set of parameters to be the reference value. Now keeping $q$ fixed,  we start with a typically high value of the discretization size $h=h_0$, we repeatedly halve it and plot the relative accuracy obtained against $h=h_0 2^{-m}$ where $m=0,1, \cdots,4$.  Relative accuracy was obtained by comparing the result obtained for a given $h$ against the reference result discussed above. The results clearly show that accuracy improves by almost 4 times with every halving of $h$.

With these choices we can replicate the results given in Johansson \& Blagouchine \cite{johansson2019computing} exactly:

\vspace{1.5mm}

$\gamma_{10^{5}} \approx
1.9919273063125410956582272431568589205211659777533113258 \\
75975525936171259272227176914320666190965225 \cdot 10^{83432},$

\vspace{1.5mm}

$\gamma_{10^{10}} \approx
7.5883621237131051948224033799125486921750410324509700470 \\
54093338492423974783927914992046654518550779 \cdot 10^{12397849705},
$

\vspace{1.5mm}

$\gamma_{10^{15}} \approx
1.8441017255847322907032695598351364885675746553315587921 \\
86085948502542608627721779023071573732022221 \cdot 10^{1452992510427658},
$

\vspace{1.5mm}

$
\gamma_{10^{100}} \approx
3.1874314187023992799974164699271166513943099108838469225 \\ 07106265983048934155937559668288022632306095 \cdot 10^e,
$

$e = 2346394292277254080949367838399091160903447689869837 \\ 3852057791115792156640521582344171254175433483694.$

\vspace{1.5mm}

Similarly for complex $v$ we obtain

\vspace{1.5mm}

$\gamma_{10^{5}}(2+3i) \approx
(1.52933142489317896667092453331813941673604063614322663 \\ 9046917471026123822028695414669890818089958104 \; + \; 7.6266053170235392288 \\
29846454534202735013368165330230700751870950104906000791927387438554979 \\ 23063058i) \cdot 10^{83440},$

\vspace{1.5mm}

$\gamma_{10^{100}}(2+3i) \approx
(2.447197253567132691871635713584630519276677767177878 \\ 733142765829147799303241971747565188937402242864 
\; + \; 1.328114485458616967 \\ 078662312208319540579816973253179511750642930437359777538176731578318799 
\\ 940692883i) \cdot 10^{e},$

$e = 23463942922772540809493678383990911609034476898698373\\852057791115792156640521582344171254175433483702$.

To calculate $\gamma_n(1)$ to 100 digits accuracy,  our method takes about 0.05 seconds for various $n$ values greater than 500.  For 1000 digits accuracy,  we take $M=3201$ and our results take about 5 seconds. These timings are at least an order of time faster than the ones reported in \cite{johansson2019computing} for $10^{100}$.  For complex $v$ the timings will be double of the corresponding timing for $v=1$ case.

\section{Asymptotic Expansion}
Applying the Saddle point approximation we obtain the following asymptotic formula for $I_n(a)$:
\begin{align}
I_n(a) = W(u)^{n+1}e^{-2  \pi  \omega} \sqrt{\frac{2 (n+1)}{\pi (1+ W(u))}} ,
\label{eq:asym}
\end{align}
where $u$ and $\omega$ are defined in \autoref{eq:owega}.  The generalized SCs, $\gamma_n(v)$, can be obtained from \autoref{eq:gammafromi}.  This is a very accurate asymptotic formula and for a given $n$,  it gives almost $\log_{10}(n)$ digits correct. This means for for $10^{10^6}$ almost an accuracy of about a million digits. For example we obtain the following results using the asymptotic formulas \\
\vspace{1.5mm}

$
\gamma_{10^{100}} \approx
3.1874314187023992799974164699271166513943099108838469225 \\ 0710626598304893415593755966828802263230609{\underline 1} \cdot 10^e,
$

$e = 2346394292277254080949367838399091160903447689869837 \\ 3852057791115792156640521582344171254175433483694,$ \\
and
\vspace{1.5mm}

$\gamma_{10^{100}}(2+3i) \approx
(2.447197253567132691871635713584630519276677767177878 \\ 73314276582914779930324197174756518893740224286{\underline 2}
\; + \; 1.328114485458616967 \\ 078662312208319540579816973253179511750642930437359777538176731578318799 
\\ 940692883i) \cdot 10^{e},$

$e = 23463942922772540809493678383990911609034476898698373\\852057791115792156640521582344171254175433483702$.

These results differ from the exact results only on the last digits which have been underlined. We note that for large $n$ other existing expansions such as that due to Knessl \& Coffey \cite{knessl2011asymptotic} and Paris \cite{paris2015asymptotic}. give similar kind of accuracy.  However, derivation provided here is considerably simpler.

\section{Conclusion}
We have shown that generalized SCs can be efficiently calculated using the DE method.  While the integration cut off $q$ can be set automatically, $M$ needs to be set manually.  It was found that M=201 is enough to obtain at least 100 digits accuracy for $n \ge 500$.   It may be possible to extend the method and set $M$ automatically.  Also it may be possible to apply modified DE method \cite{tyagi2022double} to obtain more accurate results.  We also note that cases where $|s| \ge n$ are considerably simpler to treat as the integrand is peaked around the origin.  In cases like this,  one can use the standard DE method. 

We have also derived asymptotic expansion for the generalized SCs.  Our expansion formulas are similar to the ones derived by Knessl \& Coffey and Fekih-Ahmed \cite{fekih2014new} but the derivation is simpler and based on quite different integral representations.  The asymptotic expansion formulas are highly accurate and can obtain almost an $n$ digit accuracy for the SC calculation of $\gamma_{10^n}(v)$.


\bigskip 
\bibliographystyle{unsrt} 
\bibliography{stieltjes}

\end{document}